\documentclass[a4paper,12pt]{article}
\usepackage{times}
\usepackage{amscd}
\usepackage{amsfonts}
\usepackage{amsmath}
\usepackage{amssymb}
\usepackage{amsthm}
\usepackage[all]{xy}
\def\af#1{\mathbb A^{#1}}
\def\al{\alpha}
\def\as#1{\renewcommand\arraystretch{#1}}
\def\aut{\op{Aut}}

\def\be{\bigskip}

\def\comb#1#2{\as{0.8}\left(\!\begin{array}{c}
#1\\#2
\end{array}\!\right)\as{1}}
\newcommand{\gl}{\op{GL}_2(k)}
\def\hh{\op{Hyp}}
\def\imp{\,\Longrightarrow\,}
\def\lra{\longrightarrow}
\def\mdt#1{\ \mbox{\tiny\rm(mod $#1$)}}
\def\op{\operatorname}
\def\xx{\mathcal{X}}
\def\xg{\op{Fix}_{\ga}(\xx)}
\def\yy{\mathcal{Y}}
\newcommand{\To}{\longrightarrow}

\newcommand{\cc}{{\mathcal C}}

\newcommand{\G}{\mathbb{G}_m}
\newcommand{\N}{\mathbb{N}}

\newcommand{\prl}{\mathbb{P}^{1}}

\newcommand{\md}[1]{\ \mbox{\rm(mod }{#1})}
\newcommand{\tq}{\,|\,}

\newcommand{\g}{\Gamma}

\newcommand{\ga}{\gamma}
\newcommand{\eps}{\epsilon}
\newcommand{\sg}{\sigma}
\newcommand{\la}{\lambda}

\def\fg{\op{Fix}_{\ga}}

\newcommand{\pg}{\op{PGL}_2(k)}
\def\pr#1{\mathbb P^{#1}}
\newcommand{\kb}{\overline{k}}
\def\fq{\mathbb{F}_q}
  {\end{eqnarray}%
}
\newenvironment{myeqnarray*}
{\setlength{\arraycolsep}{1.5pt}\begin{eqnarray*}}%
{\end{eqnarray*}}%

\newtheorem{theorem}{Theorem}[section]
\newtheorem{corollary}[theorem]{Corollary}

\newtheorem{remark}[theorem]{Remark}
\newtheorem{lemma}[theorem]{Lemma}
\newtheorem{proposition}[theorem]{Proposition}
\newtheorem{definition}[theorem]{Definition}

\title{Counting hyperelliptic curves}

\author{Enric Nart \thanks{partially supported by grant MTM2006-11391 from the Spanish MEC}}

\begin{document}
\maketitle

\begin{abstract}
We find a closed formula for the number $\op{hyp}(g)$ of hyperelliptic curves of genus $g$ over a finite field
$k=\fq$ of odd characteristic. These numbers $\op{hyp}(g)$ are expressed as a polynomial in $q$ with integer coefficients that depend on $g$ and the set of divisors of $q-1$ and $q+1$. As a by-product we obtain a closed formula for the number of self-dual curves of genus $g$. A hyperelliptic curve is self-dual if it is $k$-isomorphic to its own hyperelliptic twist.
\end{abstract}

\section*{Introduction} In this paper we find a closed formula for the number $\op{hyp}(g)$ of hyperelliptic curves of genus $g$ over a finite field
$k=\fq$ of odd characteristic (Theorem \ref{formula}). We use a general technique for enumerating $\pg$-orbits of rational $n$-sets of $\pr1$ that was developed in \cite{lmnx} and extended to arbitrary dimension in \cite{mn}. For $n=2g+2$, each $n$-set $S=\{t_1,\dots,t_n\}$ of $\pr1$ determines a family of hyperelliptic curves of genus $g$ whose Weierstrass points have $x$-coordinate in $S$: 
$$C_{\la,S}\colon \quad y^2=\la \prod_{t\in S,\,t\ne\infty}(x-t), \quad\la\in k^*.$$ If $S$ is a {\it rational $n$-set} (i.e. stable under the action of the Galois group $\op{Gal}(\kb/k)$), the curve $C_{\la,S}$ is defined over $k$. These curves fall generically into two different $k$-isomorphism classes, represented by a curve and its hyperelliptic twist. However, there are $n$-sets for which these curves are all $k$-isomorphic; in other words there are {\it self-dual curves} that are $k$-isomorphic to their own hyperelliptic twist. Finally, two different rational $n$-sets of $\pr1$ determine the same family of hyperelliptic curves up to $k$-isomorphism if and only if they are in the same orbit under the natural action of $\pg$. Summing up, if $\hh_g$ is the set of $k$-isomorphism classes of hyperelliptic curves over $k$ of genus $g$, and $\comb{\pr1}n(k)$ is the set of rational $n$-sets of $\pr1$ we have a well defined map
$$
\hh_g\lra \pg\backslash\comb{\pr1}{2g+2}(k),  
$$
sending a curve $C$ to the $2g+2$-set of $x$-coordinates of the Weierstrass points of a Weierstrass model of $C$.
This map is onto and the orbit of each $n$-set $S$ has either one or two preimages according to $C_{\la,S}$ being self-dual or not.

In sections 2, 3 we study when a concrete $k$-automorphism of $\pr1$ can determine a $k$-isomorphism between a curve and its hyperelliptic twist (Theorem \ref{eps}). This result provides a way of counting $\op{hyp}(g)=|\hh_g|$ by using the techniques of \cite{lmnx} and \cite{mn}, where a closed formula for the cardinality of the target set $\pg\backslash\comb{\pr1}{2g+2}(k)$ was found. As a by-product we obtain also a closed formula for the number of self-dual curves of a given genus (Theorem \ref{sd}).  \be  

\noindent{\bf Acknowledgement. }It is a pleasure to thank Amparo L\'opez for providing the key argument to prove Lemma \ref{orbit}.\be

\noindent{\bf Notations. }We fix once and for all a finite field $k=\fq$ of odd characteristic $p$
and an algebraic closure $\kb$ of $k$. We denote by $G_k$ the
Galois group $\op{Gal}(\kb/k)$ and by $\sg(x)=x^q$ the
Frobenius automorphism, which is a topological generator of $G_k$ as a profinite group. Also, $k_2$ will denote the unique quadratic extension of $k$ in $\kb$.

\section{Generalities on hyperelliptic curves}
Let $C$ be a hyperelliptic curve defined over $k$; that is,
 $C$ is a smooth, projective and geometrically irreducible curve
defined over $k$, of genus $g\ge 2$, and it admits a degree two morphism
$\pi\colon C \To \prl$, which is  also
defined over $k$.

To the non-trivial element of the Galois group of the quadratic extension of function fields, $
k(C)/\pi^*(k(\prl))$, it corresponds an
involution, $\iota\colon C\longrightarrow C$, which is called the
{\it hyperelliptic involution} of $C$. Let us recall some basic properties of $\iota$.

\begin{theorem}\label{restr}
  Let  $\pi_1 ,\pi_2 \colon C \To \prl$ be two $k$-morphisms of degree two.
  Then, there exists a unique $k$-automorphism $\gamma$ of $\prl$ such
that $\pi_2=\gamma\circ \pi_1$. In particular, the involution $\iota$ is canonical, its fixed points are the Weierstrass points of $C$, and they are the ramification points of any morphism of degree two from $C$ to $\prl$.
\end{theorem}

We denote by $\aut(C)$ the $k$-automorphism group of $C$ and by
$\aut_{\kb}(C)$ the full automorphism group of $C$.  The group of reduced geometric automorphisms of $C$ is
$$
\aut'_{\kb}(C):=\{\ga\in \aut_{\kb}(\pr1)\tq \ga(\pi(W))=\pi(W)\}, 
$$where $W$ is the set of Weierstrass points of $C$. We denote by $\aut'(C)$ the subgroup of reduced automorphisms defined over $k$.

Any automorphism $\varphi$ of $C$ fits into a commutative
  diagram:
  \begin{center}
\leavevmode
\xymatrix{ C \ar[r]^\varphi \ar[d]^\pi & C \ar[d]^\pi \\
\pr1 \ar[r]^{\ga} & \pr1
}
\end{center}
  for certain uniquely determined reduced automorphism $\ga$. The map $\ \varphi\mapsto \ga$ is a
group  homomorphism (depending on $\pi$) and we have a central exact sequence of
 groups compatible with Galois action:
  $$
  1\lra\{1,\iota\} \lra \aut_{\kb}(C) \stackrel{\pi}\lra\aut'_{\kb}(C)\lra 1.
  $$
This leads to a long exact sequence of Galois cohomology sets:
$$
1\to\{1,\iota\} \to \aut(C) \stackrel{\pi}\to\aut'(C)\stackrel{\delta}\to H^1(G_k,\{1,\iota\})\to  H^1(G_k,\aut_{\kb}(C))
$$
The $\kb/k$-twists of $C$ are parameterized by the pointed set $H^1(G_k,\aut_{\kb}(C))$
and, since $k$ is a finite field, a 1-cocycle is determined just by the choice of a $\kb$-automorphism of $C$. The twist corresponding to $\iota$ is called  the {\it hyperelliptic twist} of $C$. We say that $C$ is {\it self-dual} if it is $k$-isomorphic to its own hyperelliptc twist. 

The group $H^1(G_k,\{1,\iota\})$ is isomorphic in a natural way to $k^*/(k^*)^2$ and the curve $C$ is self-dual if and only if the two elements in this group have the same image in $H^1(G_k,\aut_{\kb}(C))$. Thus, from the above exact sequence we deduce the following criterion for self-duality:

\begin{lemma}\label{dta}
The curve $C$ is self-dual if and only if the homomorphism $\delta$ does not vanish identically.
\end{lemma}

By fixing three different points $\infty,\,0,\,1 \in \prl (k)$, we get identifications
$$k(\prl)= k(x),\qquad \op{Aut}_k(\prl)= \pg.$$
The function
field $k(C)$ is then identified to a quadratic extension of
$k(x)$ and it admits a generator $y\in k(C)$
satisfying $\ y^2= f(x)$, for some separable polynomial $f(x) \in k[x]$.

Conversely, for any separable polynomial $f(x)\in k[x]$, the Weierstrass equation
$y^2=f(x)$ determines a plane non-singular affine
curve $C^{\mbox{\tiny af}}$ defined over
$k$.
Let $C$ be the projective, smooth curve
obtained as the normalization of the projective closure of
$C^{\mbox{\tiny af}}$. The projection on the first coordinate,
$(x,y)\mapsto x$, lifts to a morphism of degree two, $\pi\colon
C\To\prl$, implicitly associated to the equation. 

The genus of the hyperelliptic curve $C$ determined by such a Weiertrass equation is $g=\lfloor (n-1)/2\rfloor$, where $n$ is  the degree of $f(x)$. 

Since $C^{\mbox{\tiny af}}$ is non-singular, we can identify the set $C^{\mbox{\tiny af}}(\kb)$
 with an open subset of $C(\kb)$; this allows us to attach
affine coordinates to all points of $C$ except for a finite set of
 points at infinity: $C(\kb)\setminus C^{\mbox{\tiny
af}}(\kb)=\pi^{-1}(\infty)$. If $n$ is odd there is a unique point at infinity, and it is a Weierstrass point defined over $k$; if $n$ is even there are two points at infinity, and they are defined over $k$ if and only if the principal coefficient of $f(x)$ is a square in $k^*$. 
  
In affine coordinates, the hyperelliptic involution is expressed as: $(x,y)^{\iota}=(x,-y)$. The Weierstrass points of $C$ are given by the roots of $f(x)$ in $\kb$, and the point at infinity if $n$ is odd. The hyperelliptic twist of $C$ has Weierstrass equation $y^2=u f(x)$, for any $u\in k^*$, $u\not\in (k^*)^2$.

\section{Classification of hyperelliptic curves up to \\$k$-isomorphism}
A {\it rational $n$-set of $\pr1$} is by definition a $k$-rational
point of the variety $\comb{\pr1}{n}$ of $n$-sets of $\pr1$. Thus,
an element $S\in\comb{\pr1}{n}(k)$ is just an unordered family
$S=\{t_1,\dots,t_n\}$ of $n$ different points of $\pr1(\kb)$, which is globally invariant under the Galois action: $S=S^{\sigma}$.

To each rational $n$-set $S$ of $\pr1$ we can attach the monic
separable polynomial $f_S(x)\in k[x]$ of degree $n$ or $n-1$ given
by:
$$
f_S(x):=\prod_{t\in S,\ t\ne\infty}(x-t).
$$
Clearly, this correspondence is multiplicative with respect to
disjoint unions:
$$
S=S_1\sqcup S_2 \imp f_S(x)= f_{S_1}(x)f_{S_2}(x).
$$

The natural action of $\pg$ on $n$-sets of $\pr1$ determines a
natural action of $\pg$ on the set of hyperelliptic curves defined
over $k$. In order to analyze this action in detail we introduce multipliers $J(\ga,S)\in k^*$ that depend in principle on the choice of a representative in $\gl$ of $\ga\in\pg$. Consider a matrix
$$
\ga=\begin{pmatrix}a&b\\c&d\end{pmatrix}\in\gl.
$$
For any $t\in\pr1(\kb)$ we can define a local multiplier
$j(\ga,t)\in \kb^*$ by
$$
j(\ga,t):=\left\{\begin{array}{ll}\det(\ga)(ct+d)^{-1}&\mbox{ if
}t\ne\infty,\,t\ne -d/c\\c&\mbox{ if }t= -d/c,\,c\ne0\\d&\mbox{ if
}t=\infty,\,c=0\\-\det(\ga)c^{-1}&\mbox{ if
}t=\infty,\,c\ne0\end{array}\right.
$$
For any rational $n$-set $S$ of $\pr1$ we define a global
multiplier $$J(\ga,S):=\prod_{t\in S}j(\ga,t)\in k^*.$$ In a
tautological way this function $J(\ga,S)$ is multiplicative with
respect to disjoint unions of rational $n$-sets. The property of $J(\ga,S)$ we are mainly interested in is the following:  

\begin{proposition}\label{pu}
For any rational $n$-set $S$ of $\pr1$ and any $\ga\in\gl$ we have
$$
(cx+d)^nf_{\ga(S)}\left(\dfrac{ax+b}{cx+d}\right)=J(\ga,S)f_S(x).
$$
\end{proposition}

\begin{proof}
All objects involved behave multiplicatively with respect to
disjoint union of $n$-sets. Thus, we need only to prove this
formula for $n=1$, where it follows by a straightforward check.
\end{proof}

We recall now the usual property of these kind of multipliers:

\begin{proposition}\label{usual}
For any rational $n$-set $S$ of $\pr1$ and any $\ga,\,\rho\in\gl$
we have
$$
J(\rho,S)J(\ga,\rho(S))=J(\ga\rho,S).
$$
\end{proposition}

\begin{proof}
It is straightforward to check that
$j(\rho,t)j(\ga,\rho(t))=j(\ga\rho,t)$ for all $t\in \pr1(\kb)$.
\end{proof}

The following facts are an immediate consequence of this
property:

\begin{corollary}\label{conj}
\begin{enumerate}
\item If $\ga(S)=S$ we have $J(\rho\ga\rho^{-1},\rho(S))=J(\ga,S)$
for all $\rho\in\gl$. \item Let $\gl_S$ be the isotropy subgroup of $S$: $$\gl_S:=\{\ga\in\gl\tq \ga(S)=S\}.$$
Then, $J(\ ,S)\colon \gl_S\lra k^*$ is a group homomorphism.
\end{enumerate}
\end{corollary}

From now on we assume that $n$ is an even integer, $n>4$. To every $\la\in k^*$, $S\in
\comb{\pr1}{n}(k)$, we can attach the hyperelliptic
curve $C_{\la,S}$ determined by the Weierstrass equation $y^2=\la
f_S(x)$.

For any $\mu\in k^*$ the morphism $(x,y)\mapsto(x,\mu y)$ sets a
$k$-isomorphism between  $C_{\la,S}$ and $C_{\la\mu^2,S}$. Thus, if we let the pairs $(\la,S)$ run on the set 
$$
(\la,S)\in\xx_n:=(k^*/(k^*)^2)\times \comb{\pr1}n(k),
$$
the curves $C_{\la,S}$ contain representatives of all
$k$-isomorphism classes of hyperelliptic curves. We shall abuse of notation and denote by the same symbol $C_{\la,S}$ a concrete curve when $\la\in k^*$ and a class of curves up to $k$-isomorphism when $\la\in k^*/(k^*)^2$.

For any $\mu\in k^*$, we have clearly $J(\mu\ga,S)=\mu^nJ(\ga,S)$.
Since $n$ is even, the class of $J(\ga,S)$ modulo $(k^*)^2$
depends only on the image of $\ga$ in $\pg$. Thus, the above properties of $J(\ga,S)$ extend in an obvious way to elements of $\pg$ as long as we are interested only in the value of $J(\ga,S)$ modulo squares. In particular, the following action of $\pg$ on the set $\xx_n$ is well defined:
$$
\ga(\la,S):=(\la J(\ga,S),\ga(S)).
$$

\begin{theorem}\label{red}
Two hyperelliptic curves $C_{\la,S}$, $C_{\mu,T}$ are
$k$-isomorphic if and only if there exists $\ga\in\pg$ such that
$(\mu,T)=\ga(\la,S)$. In particular, the map $(\la,S)\mapsto
C_{\la,S}$ induces a 1-1 correspondence
$$
\pg\backslash\xx_{2g+2} \lra \hh_g,
$$
where $\hh_g$ is the set of $k$-isomorphism classes of hyperelliptic curves over $k$ of genus $g$.
\end{theorem}

\begin{proof} 
By Proposition \ref{pu} the following map is a $k$-isomorphism betwen $C_{\la,S}$ and
$C_{\la J(\ga,S),\ga(S)}$:
$$ (x,y)\mapsto
\left(\dfrac{ax+b}{cx+d},\dfrac{J(\ga,S)y}{(cx+d)^{n/2}}\right).
$$

Suppose now that the curves $C_{\la,S}$,
$C_{\mu,T}$ are $k$-isomorphic. By Theorem \ref{restr} any $k$-isomorphism
from $C_{\la,S}$ to $C_{\mu,T}$ has the shape
$$
(x,y)\mapsto (\ga(x),\psi(x,y)),
$$
for a uniquely determined $\ga\in\pg$ and certain rational
function $\psi(x,y)$ defined over $k$. Let us choose any
representative of $\ga$ in $\gl$, which we still denote by $\ga$.
By Proposition \ref{pu}, we nave necessarily $\psi(x,y)=u
J(\ga,S)y/(cx+d)^{n/2}$ for certain undetermined constant $u\in
k^*$; thus, $(\mu,T)=\ga(\la,S)$.
\end{proof}

\section{Computation of $J(\ga,S)$ modulo squares}
Our aim is the computation of $|\hh_g|$. By Theorem \ref{red} this is equivalent to the computation 
of $|\pg\backslash\xx_{2g+2}|$. In order to compute the number of these orbits it is crucial to
compute $J(\ga,S)$ modulo $(k^*)^2$ when $\ga(S)=S$. This will be achieved in Theorem \ref{eps}, where we show that the value of $J(\ga,S)$ modulo squares depends only on the order of $\ga$ as an element of $\pg$. As a first glimpse of this fact we check first what happens when $\ga$ has odd order. 

For any rational $n$-set $S$ of $\pr1$ we denote by $$\pg_S:=\{\ga\in\pg\tq \ga(S)=S\}$$ the isotropy subgroup of $S$. 
By Corollary \ref{conj} we have a homomorphism
$$
J(\ ,S)\colon \pg_S \lra  k^*/(k^*)^2,
$$  
which is a reinterpretation of the homomorphism $\delta$ of Lemma \ref{dta}, after natural identifications $\pg_S=\aut'(C)$, $H^1(G_k,\{1,\iota\})=k^*/(k^*)^2$. In particular, the curve $C_{\la,S}$ is self-dual if and only if there exists $\ga\in \pg_S$ such that $J(\ga,S)$ is not a square. Since $k^*/(k^*)^2$ is a $2$-torsion group this homomorphism vanishes on elements of odd order.

\begin{corollary}\label{modd}
The element  $J(\ga,S)$ is a square in $k^*$ for all $\ga\in\pg_S$ of odd order.
\end{corollary}

By Corollary \ref{conj}, if $\ga\in\pg_S\,$ it is sufficient to compute $J(\ga,S)$ modulo squares for a system of representatives of conjugacy classes
of $\pg$. Let us recall how these representatives can be chosen and the possible values of the order of $\ga$ in each conjugacy class.

\begin{remark}\label{rmk}
 Let $\ga\in\pg$, $\ga\ne1$, and let $m$ be the order of $\ga$. Let $\fg$ denote the set of fixed points of
$\ga$ in $\pr1(\kb)$. 

There are three possibilities for the conjugacy class of $\ga$:
\begin{enumerate}
\item The automorphism $\ga$ is homothetic (conjugate to $t\mapsto \la t$, $\la
\in k^*$). Then $\fg$ consists of two points in $\pr1(k)$; in this
case $m$ is the order of $\la$ in $k^*$, which is a divisor of
$q-1$. \item The automorphism $\ga$ is conjugate to the translation, $t\lra t+1$.
Then $\fg$ consists of a single point in $\pr1(k)$; in this case
$m=p$. \item The automorphism  $\ga$ is potentially homothetic; i.e. $\ga$ is conjugate to the class in $\pg$ of a matrix
$\begin{pmatrix}0&1\\c&d\end{pmatrix}\in\gl$ with eigenvalues
$\al,\,\al^{\sg}$ in $k_2\setminus k$. Then $\fg$ consists of two quadratic
conjugate points in $\pr1(k_2)$; in this case $m$ is the least
positive integer such that $\al^m\in k$, and it is a divisor of
$q+1$.
\end{enumerate}
\end{remark}

Let us introduce
a special notation for the class of $J(\ga,S)$ modulo squares.

\begin{definition}
Denote by $\eps$ the map
$$\eps\colon\pg\times \comb{\pr1}n (k)\stackrel{J}{\lra} k^* \lra k^*/(k^*)^2 \lra \{\pm1\},
$$
where the last map is the unique non-trivial group homomorphism
between these two groups of order two.
\end{definition}

\begin{theorem}\label{eps}
Let $S$ be a rational $n$-set of $\pr1$ and let $\ga\in\pg_S$ of order $m$. Then, if $n$ is even: 
$$\as{1.2}
\eps(\ga,S)=\left\{ \begin{array}{ll} (-1)^{(q-1)/m},&\mbox{ if
$\ga(t)=\la t,\  \la\in k^*,$  and }\infty\in S,\\
1,&\mbox{ if
$\ga(t)=\la t,\  \la\in k^*,$  and }\infty\not\in S,\\
1,&\mbox{ if }
\ga(t)=t+1,\\(-1)^{(q+1)/m}(-1)^{(n-2)/m},&\mbox{ if
$\ga$ potent. homothetic and }\fg\subseteq S,
\\
(-1)^{n/m},&\mbox{ if
$\ga$ potent. homothetic and }\fg\not\subseteq S.
\end{array}\right.
$$
\end{theorem}

In order to prove this theorem we need a couple of lemmas.

\begin{lemma}\label{orbit}  Let $\ga$ be a potentially homothetic element of $\pg$ of order $m$, with representative $\begin{pmatrix}0&1\\c&d\end{pmatrix}$ in $\gl$. Let $\al,\,\al^{\sg}\in k_2\setminus k$ be the eigenvalues of this matrix.
Let $t\in \pr1(\kb)$ be such that $\ga(t)\ne t$ and denote by $\op{O}_{\ga}(t)$ the orbit of $t$ under the action of the cyclic group generated by $\ga$. Then, $J(\ga,\op{O}_{\ga}(t))=\al^m$. 
\end{lemma}

\begin{proof}These orbits $\op{O}_{\ga}(t)=\{\ga^i(t)\tq i\in\N\}$ have either one element, if $t$ is a fixed point of $\ga$, or exactly $m$ elements (and not less) otherwise (\cite[Lem. 2.3]{lmnx} or \cite[Lem. 4.2]{mn}).
Note that the three points $\infty,\,0,\,-d/c$ (two points if $m=2$ because then $d=0$) are in the same orbit under the action of the cyclic group 
generated by $\ga$. Suppose first that $\infty\not\in  O_{\ga}(t)$; then, for all $s$ in the orbit we have $j(\ga,s)=\det(\ga)/(cs+d)=-c\ga(s)$ and
$$
J(\ga,O_{\ga}(t))=(-c)^m\prod_{s\in O_{\ga}(t)}\ga(s)=(-c)^m\prod_{s\in O_{\ga}(t)}s.
$$
In order to compute the last product  we observe that the second row of $\ga^i$ coincides with the first row of $\ga^{i+1}$: 
$$
\begin{pmatrix}0&1\\c&d\end{pmatrix}
\begin{pmatrix}r&s\\u&v\end{pmatrix}=
\begin{pmatrix}u&v\\\star&\star\end{pmatrix},
$$
and $\ga^m=\begin{pmatrix}\al^m&0\\0&\al^m\end{pmatrix}$. Thus in the product of all elements in the orbit every denominator cancels with the numerator of the following term:
$$
\prod_{s\in O_{\ga}(t)}s=t\cdot \dfrac1{ct+d}\cdot \dfrac{ct+d}{ut+v}\cdot \dfrac{ut+v}{\star\, t+\star}\cdot\cdots\cdot
\dfrac{\star\, t+\star}{\al^m t}=\al^{-m},$$the last denominator being equal to the numerator of $\ga^m(t)=\al^mt/\al^m$. Since $\al\al^{\sg}=-c$, we have $(-c)^m=\al^{2m}$, so that $J(\ga,O_{\ga}(t))=\al^m$, as claimed in the statement of the lemma. 

Suppose now that our orbit $\op{O}_{\ga}(t)$ contains $\infty,\,0,\,-d/c$, and $m>2$. Then,
$$
\prod_{s\in \op{O}_{\ga}(t),\,s\ne\infty,\,-d/c}j(\ga,s)=\prod_{s\in \op{O}_{\ga}(t),\,s\ne\infty,\,-d/c}(-c\ga(s))=(-c)^{m-2}\prod_{i=1}^{m-2}\ga^i(0).
$$  Since $ \ga^{m-1}=\al^m\ga^{-1}=\al^m\begin{pmatrix}-dc^{-1}&c^{-1}\\1&0\end{pmatrix}$, the last product can be computed as above:
$$
\prod_{i=1}^{m-2}\ga(0)=\dfrac1{c\cdot 0+d}\cdot \dfrac{c\cdot 0+d}{u\cdot 0+v}\cdot \dfrac{u\cdot 0+v}{\star\cdot 0 +\star}\cdot\cdots\cdot
\dfrac{\star \cdot 0+\star}{-\al^m dc^{-1}\cdot 0+\al^{m}c^{-1}}=c\al^{-m},
$$the last denominator given by the first row of  $ \ga^{m-1}$. Since $j(\ga,\infty)=1$, and $j(\ga,-d/c)=c$, we get 
$$J(\ga,O_{\ga}(0))=(-c)^{m-2}c\al^{-m}j(\ga,\infty)j(\ga,-d/c)=(-c)^m\al^{-m}=\al^m,$$
as claimed. In the case $m=2$ we have $d=0$ and we get the same result $$J(\ga,O_{\ga}(\infty))=j(\ga,\infty)j(\ga,0)=1\cdot c=\al^2.$$
\end{proof}

\begin{lemma}\label{norm}  Let  $\al$ be an element of $k_2\setminus k$ and  let $m>1$ be the order of $\al$ in the cyclic group $k_2^*/k^*$; in other words, $m$ is the least positive integer such that $\al^m$ belongs to $k^*$. Then,
\begin{enumerate}
\item The norm  $\op{N}_{k_2/k}(\al)$ is a square in $k^*$ if and only if $(q+1)/m$ is even.
\item If $m$ is even, the element $\al^m$ is not a square in $k^*$. 
\end{enumerate}
\end{lemma}

\begin{proof}
Note that $\op{N}_{k_2/k}(\al)=\al^{q+1}$ and $(\al^{q+1})^{(q-1)/2}=\al^{(q^2-1)/2}$; hence, $\op{N}_{k_2/k}(\al)$ is a square in $k^*$ if and only if $\al$ is a square in $k_2^*$. On the other hand, let $\al=\zeta^u$, where $\zeta$ is a generator of the cyclic group $k_2^*$. Since $\al^m=\zeta^{mu}\in k$ and the class of $\zeta$ modulo $k^*$ is a generator of $k_ 2^*/k^*$ we have necessarily $mu\equiv 0 \md{q+1}$ and $m$ is minimum with this property; that is, $m=(q+1)/\gcd(u,q+1)$. Hence, $(q+1)/m=\gcd(u,q+1)$ is even if and only if $u$ is even, and this is equivalent to the fact that $\al$ is a square in $k_2^*$.  This proves the first item.

We have $um=(q+1)v$, with $v=u/\gcd(u,q+1)$. Now,
$$
(\al^m)^{(q-1)/2}=\zeta^{um(q-1)/2}=\left(\zeta^{(q^2-1)/2}\right)^v=(-1)^v,
$$
and if $m$ is even, we have necessarily $v_2(u)<v_2(q+1)$ and this implies that $v$ is odd. This proves item 2.
\end{proof}

\noindent{\bf Proof of Theorem \ref{eps}. }
If $m$ is odd the theorem claims that $\eps(\ga,S)=1$, in agreement with Corollary \ref{modd}. From now on we suppose that $m$ is even.
 
Suppose $\ga(t)=\la t$ for some $\la\in k^*$. By choosing the matrix $\begin{pmatrix}\la&0\\0&1\end{pmatrix}$ as a representative of $\ga$ we get $j(\ga,t)=\la$ if $t\ne\infty$, and $j(\ga,\infty)=1$. Thus, $$J(\ga,S)=\left\{\begin{array}{ll}
\la^n,&\mbox{ if }\infty\not \in S,\\\la^{n-1},&\mbox{ if }\infty \in S.\end{array}\right.$$    
Hence, $\eps(\ga,S)=1$ if $\infty\not \in S$. On the other hand, if $\infty\in S$ we get $\eps(\ga,S)=1$ if and only if $\la$ is a square in $k^*$; this is equivalent to $m$ being a divisor of $(q-1)/2$, which is in turn equivalent to $(q-1)/m$ being even.

Suppose $\ga(t)=t+1$. We can choose the matrix $\begin{pmatrix}1&1\\0&1\end{pmatrix}$ as a representative of $\ga$, leading to $j(\ga,t)=1$ for all $t\in\pr1(\kb)$. Hence $\eps(\ga,S)=1$.

Suppose $\ga$ potentially homothetic with representative $\begin{pmatrix}0&1\\c&d\end{pmatrix}$ in $\gl$ and let $\al,\,\al^{\sg}\in k_2\setminus k$ be the eigenvalues of this matrix. Since $S$ is invariant under the action of $\ga$, it is a disjoint union of orbits of the cyclic group generated by $\ga$. We can compute $J(\ga,S)$ as the product of the different $J(\ga,O_{\ga}(t))$. 

If $t=1/(ct+d)$ is a fixed point of $\ga$, we have $j(\ga,t)=-c/(ct+d)=-ct$. If $S$ contains a fixed point of $\ga$ then necessarily $\fg\subseteq S$, because $S$ is invariant under the action of Frobenius; in this case the  multiplier $J(\ga,S)$ has a factor $j(\ga,\fg)=(-ct)(-ct^{\sg})=c^2tt^{\sg}=-c=\op{N}_{k_2/k}(\al)$, because the minimal polynomial of $t$ is $x^2+dc^{-1}-c^{-1}$. 

Therefore, if $\fg\not\subseteq S$, $S$ is the union of $n/m$ orbits of cardinality $m$ and Lemma \ref{orbit} shows that
$J(\ga,S)=(\al^m)^{n/m}$; by Lemma  \ref{norm} this is a square in $k^*$ if and only if $n/m$ is even. On the other hand, if $\fg\subseteq S$ then $S$ is the union of $\fg$ and $(n-2)/m$ orbits of cardinality $m$ and Lemma \ref{orbit} shows that $J(\ga,S)=\op{N}_{k_2/k}(\al)(\al^m)^{(n-2)/m}$; by Lemma  \ref{norm} this is a square in $k^*$ if and only if $(-1)^{(q+1)/m}(-1)^{(n-2)/m}=1$.\hfill{$\Box$}\bigskip

\section{Counting $\pg$-orbits of $\xx_{2g+2}$}
Let $\g$ be a finite group acting on a finite set $\xx$. The number of orbits of this action can be counted as the average number of fixed points:
\begin{equation}\label{cf} |\g\backslash \xx|=\frac
1{\vert\g\vert}\sum_{\ga\in\g}\vert \xg\vert=
\sum_{\ga\in\cc}\frac{\vert \xg\vert}{\vert \g_{\ga}\vert},
\end{equation} where  $\cc$ is a set of representatives of conjugacy classes of elements of
$\g$ and
 $$
 \xg:=\{x\in \xx \tq \ga(x)=x\}, \quad
 \g_{\ga}:=\{\rho\in \g \tq \rho\ga\rho^{-1}=\ga\}.
 $$

In this section we apply this formula to compute the number $\op{hyp}(g)$ of orbits of $\xx:=\xx_{2g+2}$ under the action of $\g:=\pg$. By Theorem \ref{red} this is the number of $k$-isomorphy classes of hyperelliptic curves of genus $g$. In \cite{lmnx} similar ideas were applied to count the total number of $\pg$-orbits of rational $n$-sets of $\pr1$. In \cite{mn} a general theory is developed to deal with similar problems in arbitrary dimension. For commodity of the reader we sum up the results we are going to use of these two papers.

The conjugacy classes of $\pg$ are divided into four {\it types} A, B, C, D, according to $\ga$ being respectively potentially homothetic, conjugated to a translation, homothetic or the identity (cf. Remark \ref{rmk}). The types B and D contain a single conjugacy class, represented by $\ga(t)=t+1$ (type B) and $\ga(t)=t$ (type D). The conjugacy classes of types A and C are divided into {\it subtypes}, according to the order $m$ of $\ga$; thus, the subtypes of type A are parameterized by divisors $m>1$ of $q+1$, and the subtypes of type C are parameterized by divisors $m>1$ of $q-1$. The computation of $\op{hyp}(g)$ given by (\ref{cf}) can be splitted into the sum of four terms $\op{hyp}(g)=h_A(g)+h_B(g)+h_C(g)+h_D(g)$, each term taking care of the contribution of all conjugacy classes in a concrete type. By Theorem \ref{eps} the cardinality $|\xg|$ will depend only on the subtype of $\ga$ (cf. (\ref{tres}) below), and the same happens with $|\g_{\ga}|$; hence, we can group together all $\ga$ in the same subtype and after an explicit computation of $|\g_{\ga}|$ and the number of $\ga$ in each subtype (given respectively in \cite[Prop. 2.3, Lem. 2.4]{ln}) we can express each of the above partial terms as:
\begin{equation}\label{abcd}\as{2.4}
\begin{array}{ll}
h_A(g)=\sum_{1<m|(q+1)}\dfrac{\varphi(m)|\xg|}{2(q+1)},\quad h_B(g)=\dfrac{|\xg|}q,\\
h_C(g)=\sum_{1<m|(q-1)}\dfrac{\varphi(m)|\xg|}{2(q-1)},\quad h_D(g)=\dfrac{|\xg|}{q(q-1)(q+1)}, 
\end{array}
\end{equation}
where $\ga$ is an arbitrary choice of an element in each subtype.

For the computation of $|\xg|$ we need to know the total number of rational $n$-sets of several quasiprojective subvarieties of $\pr1$. For any quasiprojective variety $V$ defined over $k$ we denote by $$a_V(n):=\left|\comb Vn(k)\right|$$ the number of rational $n$-sets of $V$. In \cite[Thm. 1.2]{mn} the generating function of these numbers in expressed in terms of the zeta function of $V$ and it is straightforward to derive from this result explicit formulas for these numbers $a_V(n)$. For our purposes we need these formulas for the cases $V=\pr1,\,\af1,\,\G$ and $\pr1_0$, where $\pr1_0$ will denote in the sequel the subvariety $\pr1\setminus\{t,t^{\sg}\}$, being $t$ any point in $\pr1(k_2)\setminus\pr1(k)$. Actually, in our formula for $\op{hyp}(g)=|\hh_g|$ there will appear certain normalizations $A_i(n)$ of these numbers. The following lemma collects all the formulas we need, extracted from \cite[Lem. 2.1]{lmnx}.  

\begin{lemma}\label{atilla}
$$
a_{\pr{1}}(n)=q^n-q^{n-2},\quad \mbox{ if } \,n\ge3.
$$

$$
A_1(n):=\dfrac{a_{\af1}(n)}q=\left\{\begin{array}{ll}
1,&\mbox{ if }n=1,\\
q^{n-1}-q^{n-2},&\mbox{ if }n\ge2.\end{array}\right.
$$
   
$$
A_2(n):=\dfrac{a_{\G}(n)}{q-1}=\dfrac{q^n-(-1)^n}{q+1},\quad
\forall n\ge 1.
$$

$$
A_0(n):=\dfrac{a_{\pr1_0}(n)}{q+1}=
\dfrac{q^{n+1}-q^n-(-1)^{\lceil n/2\rceil}q+(-1)^{\lceil(n-1)/2\rceil}}{q^2+1}.
$$
\end{lemma}

Finally, we need \cite[Cor. 2.4]{mn} that computes the number of rational $n$-sets of $V$ that are invariant under the action of $\ga$.

\begin{lemma}\label{quot}
Let $\ga\in\pg$ of order $m$ and let $V\subseteq \pr1$ be a quasiprojective subvariety which is invariant by $\ga$ and contains no fixed points of $\ga$. Then,
$$
\left|\fg\left(\comb V{nm}(k)\right)\right|=a_{V/\ga}(n)=a_V(n).
$$ 
\end{lemma}

With these results and Theorem \ref{eps} it is easy to count $\op{hyp}(g)$. In general, a pair $(\la,S)\in\xx$ is invariant by $\ga$ if and only if $\ga(S)=S$ and $\eps(\ga,S)=1$. Thus,
\begin{equation}\label{tres}
\left|\fg(\xx)\right|=2\left|\left\{S\in \fg\left(\comb{\pr1}{2g+2}(k)\right)\tq \eps(\ga,S)=1\right\}\right|.
\end{equation}

We consider separatedly the cases of $\ga$-invariant $2g+2$-sets not containing fixed points of $\ga$, $\ga$-invariant $2g+1$-sets together with one rational fixed point of $\ga$ and $\ga$-invariant $2g$-sets with two fixed points. By Lema \ref{quot}, Theorem \ref{eps} and (\ref{abcd}) we get

\begin{small}
\begin{multline*}
h_A(g)=\sum_{1<m|q+1}\dfrac{\varphi(m)}{q+1}\left(\left[a_{\pr1_0}\left(\dfrac {2g+2}m\right)\right]_{\frac {2g+2}m \mbox{\tiny \,even}}+
\left[a_{\pr1_0}\left(\dfrac {2g}m\right)\right]_{\frac {2g}m \equiv \frac {q+1}m\mdt2}\right)=\\=\sum_{1<m|q+1}\varphi(m)\left(\left[A_0\left(\dfrac {2g+2}m\right)\right]_{\frac {2g+2}m \mbox{\tiny \,even}}+
\left[A_0\left(\dfrac {2g}m\right)\right]_{\frac {2g}m \equiv \frac {q+1}m\mdt2}\right).
\end{multline*}
 
$$
h_B(g)=\dfrac 2q \left(a_{\af1}\left(\frac {2g+2}p\right)+a_{\af1}\left(\frac{2g+1}p\right)\right)=2A_1\left(\frac {2g+2}p\right)+2A_1\left(\frac{2g+1}p\right).
$$

\begin{multline*}
h_C(g)=\sum_{1<m|q-1}\dfrac{\varphi(m)}{q-1}\left(a_{\G}\left(\frac{2g+2}m\right)+2a_{\G}\left(\frac{2g+1}m\right)+\right.\\
\left.+\left[a_{\G}\left(\frac{2g}m\right)\right]_{\frac {q-1}m \mbox{\tiny\, even}}\right)=\\=\sum_{1<m|q-1}\varphi(m)\left(A_2\left(\frac{2g+2}m\right)+2A_2\left(\frac{2g+1}m\right)+
\left[A_2\left(\frac{2g}m\right)\right]_{\frac {q-1}m \mbox{\tiny\, even}}\right).
\end{multline*}

$$
h_D(g)=\dfrac{2a_{\pr1}(n)}{q(q-1)(q+1)}=2q^{2g-1}.
$$
\end{small}

By convention, in these formulas we consider $a_V(x)=0=A_i(x)$ if $x$ is not a positive integer. Also, a term $[x]_{\mbox{\tiny condition}}$ in a formula means: add $x$ if the ``condition" is satisfied.

By using the explicit formulas for $A_i(n)$ given in Lemma \ref{atilla} we obtain a closed formula for $\op{hyp}(g)$ as a polynomial in $q$ with integer coefficients that depend on the set of divisors of $q-1$ and $q+1$. This is more clearly seen if we rewrite our computation of $\op{hyp}(g)$ in a way that is more suitable for an effective computation when $g$ is given and we want to deal with a generic value of $q$.

\begin{theorem}\label{formula}
\begin{multline*}
\op{hyp}(g)=2q^{2g-1}+\sum_{1<m|2g+2}\left(\varphi(m)\left[A_0\left(\dfrac {2g+2}m\right)\right]_{m|q+1,\,\frac {2g+2}m \mbox{\tiny \,even}}+\right.\\\left.+\varphi(m)\left[A_2\left(\frac{2g+2}m\right)\right]_{m|q-1}+2\left[A_1\left(\frac{2g+2}m\right)\right]_{m=p}\right)+\\
+\sum_{1<m|2g+1}\left(2\varphi(m)\left[A_2\left(\frac{2g+1}m\right)\right]_{m|q-1}+2\left[A_1\left(\frac{2g+1}m\right)\right]_{m=p}\right)+\\+
\sum_{1<m|2g}\varphi(m)\left(\left[A_0\left(\dfrac {2g}m\right)\right]_{m|q+1,\,\frac {2g}m \equiv \frac {q+1}m\mdt2}+\left[A_2\left(\frac{2g}m\right)\right]_{m|q-1,\,\frac {q-1}m \mbox{\tiny\, even}}\right).\end{multline*}
\end{theorem}

\begin{center}\as{1.4}
\begin{table}
\caption{\small Number of hyperelliptic curves of genus $g$ up to $k$-isomorphism}\as{1.2}
\begin{tabular}{|c|l|}\hline
\mbox{\footnotesize$g$}&\qquad\qquad\qquad\qquad\qquad\ \mbox{\footnotesize$\op{hyp}(g)=|\hh_g|$}\\\hline
\mbox{\footnotesize$2$}&\mbox{\footnotesize$2q^3+q^2+2q-2+[2]_{3|q-1}+[8]_{5|q-1}+[2]_{p=5}+[2]_{q\equiv1,\,3\mdt8}$}\\\hline
\mbox{\footnotesize$3$}&\negmedspace\negmedspace\mbox{\footnotesize$\begin{array}{l}2q^5+2q^3-2-2[q^2-q]_{4|q+1}+2[q-1]_{p>3}+[4]_{8|q-1}+[12]_{7|q-1}+[2]_{p=7}+\\\qquad\qquad\qquad\qquad\qquad\qquad+[2]_{q\equiv1,\,5\mdt{12}}\end{array}$}\\\hline
\mbox{\footnotesize$4$}&\negmedspace\negmedspace\mbox{\footnotesize$\begin{array}{l}
2q^7+q^4+4[q^2-q+1]_{3|q-1}+2[q^2-q]_{p=3}-2[q-1]_{4|q+1}+4[q-1]_{q\equiv\pm1\mdt5}+\\
\qquad\quad+2[q-1]_{p=5}+2[q-1]_{q\equiv\pm1\mdt8}+[4]_{5|q-1}+[12]_{9|q-1}+[4]_{q\equiv1,\,7\mdt{16}}\end{array}$}\\\hline
\mbox{\footnotesize$5$}&\negmedspace\negmedspace\mbox{\footnotesize$\begin{array}{l}
2q^9+2q^5+2-2[q^4-q^3+q^2-q+2]_{4|q+1}+4[q-1]_{3|q-1}+4[q-1]_{q\equiv\pm1\mdt5}+\\
\qquad\qquad\qquad\qquad\qquad\qquad+[4]_{12|q-1}+[20]_{11|q-1}+[2]_{p=11}+[4]_{q\equiv1,\,9\mdt{20}}\end{array}$}\\\hline
\mbox{\footnotesize$6$}&\negmedspace\negmedspace\mbox{\footnotesize$\begin{array}{l}
2q^{11}+q^6-2[q^3-q^2]_{4|q+1}+2[q^3-q^2+q-1]_{3|q-1}+2[q^3-q^2-q+1]_{3|q+1}+\\
\qquad+2[q^2-q+1]_{8|q-1}+2[q^2-q-1]_{8|q-3}+2[q-1]_{q\equiv\pm1\mdt{12}}+2[q-1]_{p=7}+\\\qquad\qquad+6[q-1]_{q\equiv\pm1\mdt7}+[6]_{7|q-1}+[2]_{p=13}+[24]_{13|q-1}+[4]_{q\equiv1,\,11\mdt{24}}\end{array}$}\\\hline
\mbox{\footnotesize$7$}&\negmedspace\negmedspace\mbox{\footnotesize$\begin{array}{l}
2q^{13}+2q^7-2q^5+4q^3-2q^2-2-2[q^6-q^5+q^2+q-2]_{4|q+1}+2[q^4-q^3]_{p=3}+\\
\quad\ +4[q^4-q^3+q^2-q+1]_{3|q-1}+2[q^2-q]_{p=5}+8[q^2-q+1]_{5|q-1}+[8]_{16|q-1}+\\\qquad\qquad+4[q-1]_{q\equiv\pm1\mdt8}+6[q-1]_{q\equiv\pm1\mdt7}+[16]_{15|q-1}+[6]_{q\equiv1,\,13\mdt{28}}\end{array}$}\\\hline
\mbox{\footnotesize$8$}&\negmedspace\negmedspace\mbox{\footnotesize$\begin{array}{l}
2q^{15}+q^8+2[q^5-q^4]_{4|q-1}-2[q-1]_{4|q+1}+2[q^3]_{3|q-1}-2[q^3-q^2-q+1]_{3|q+1}+\\
\qquad+2[q^3-q^2+q-1]_{8|q-1}+2[q^3-q^2-q+1]_{8|q+1}+6[q-1]_{q\equiv\pm1\mdt9}+\\
\qquad\qquad+4[q-1]_{q\equiv\pm1\mdt{16}}+[6]_{9|q-1}+[32]_{17|q-1}+[2]_{p=17}+[8]_{q\equiv1,\,15\mdt{32}}\end{array}$}\\\hline
\mbox{\footnotesize$9$}&\negmedspace\negmedspace\mbox{\footnotesize$\begin{array}{l}
2q^{17}+2q^9-2q^8+2q^5-2q^4+2q-2+2[q^5-q^4+q-1]_{p>3}2+\\
\quad+[q^8-q^7+2q^4-2q^3+q^2-q+2]_{4|q-1}+2[q^3-q^2]_{3|q-1}-2[q^3-q^2]_{3|q+1}+\\\qquad+2[q^3-q^2]_{p=5}+4[q^3-q^2]_{q\equiv\pm1\mdt5}+
2[q^2-q]_{q\equiv1,\,5\mdt{12}}+\\\qquad\qquad+8[q-1]_{5|q-1}+6[q-1]_{q\equiv\pm1\mdt9}+[8]_{20|q-1}+[36]_{19|q-1}+[2]_{p=19}+\\
\qquad\qquad\qquad\qquad+[2]_{12|q-1}-[2]_{12|q-5}
+[6]_{q\equiv1,\,17\mdt{36}}\end{array}$}\\\hline
\mbox{\footnotesize$10$}&\negmedspace\negmedspace\mbox{\footnotesize$\begin{array}{l}
2q^{19}+q^{10}-2[q^7-q^6+q^3-q^2]_{4|q+1}+4[q^6-q^5+q^4-q^3+q^2-q+1]_{3|q-1}+\\\quad+2[q^6-q^5]_{p=3}+2[q^4-q^3+q^2-q+1]_{8|q-1}+2[q^4-q^3-q^2+q+1]_{8|q-3}+\\\qquad\quad+4[q^3-q^2+q-1]_{5|q-1}+4[q^3-q^2-q+1]_{5|q+1}+12[q^2-q+1]_{7|q-1}+\\
\qquad\qquad\qquad+2[q^2-q]_{p=7}+10[q-1]_{11|q+1}+[10q]_{11|q-1}+2[q-1]_{p=11}+\\\qquad\qquad\qquad\qquad+4[q-1]_{q\equiv\pm1\mdt{20}}+[24]_{21|q-1}+[8]_{q\equiv1,\,19\mdt{40}}\end{array}$}\\\hline
\end{tabular}
\end{table}
\end{center}

We display in Table 1 the value of $\op{hyp}(g)$ for $2\le g\le 10$. For $g=2$, Theorem \ref{formula} corrects
the formula for $\op{hyp}(2)$ given in \cite[Thm. 22]{car}, which has a right generic term but a wrong term depending on the class of $q$ modulo $120$.
 
On the other hand, just by considering the contribution of $m=2$ in the formula of Theorem  \ref{formula} we get a fairly good general approximation to $\op{hyp}(g)$:

\begin{corollary}
$$\as{1.2}
\begin{array}{ll}
\op{hyp}(g)=2q^{2g-1}+q^g+O(q^{\lfloor(2g-1)/3\rfloor}),&\qquad \mbox{ if $g$ even, }g\ge 8,\ 4|q-1, \\
\op{hyp}(g)=2q^{2g-1}+q^g-2q^{g-3}+O(q^{g-4}),&\qquad \mbox{ if $g$ even, }g\ge 10,\ 4|q+1, \\
\op{hyp}(g)=2q^{2g-1}+2q^g-2q^{g-2}+O(q^{g-4}),&\qquad \mbox{ if $g$ odd,  }g\ge 9,\ 4|q-1,\\  
\op{hyp}(g)=2q^{2g-1}+2q^g-2q^{g-1}+O(q^{g-4}),&\qquad \mbox{ if $g$ odd, }g\ge 9,\ 4|q+1.
\end{array}
$$
\end{corollary}

\begin{proof}
Apart from the generic term $2q^{2g-1}$, the contribution of $m=2$ in the general formula is:
$$
[A_0(g+1)]_{g\mbox{\tiny\, odd}}+A_2(g+1)+[A_0(g)]_{g\equiv(q+1)/2\mdt2}+[A_2(g)]_{4|q-1}.
$$
The result follows by applying Lemma \ref{atilla}, having in mind that the highest power of $q$ arising from the other terms is $q^{\lfloor(2g-1)/3\rfloor}$, corresponding to $m=3$.
\end{proof}

\section{Self-dual curves}
In this section we find a closed formula for the number $\op{sd}(g)$ of self-dual hyperelliptic curves of genus $g$, up to $k$-isomorphism.  As mentioned in the Introduction,
$$
\op{hyp}(g)=\sum_{S\in\yy}\delta_S,\qquad \yy:=\pg\backslash\comb{\pr1}{2g+2}(k),
$$
where $\delta_S=1$ if the curve $y^2=f_S(x)$ is self-dual and $\delta_S=2$ otherwise. Thus, 
$$
\op{sd}(g)=2\left|\yy\right|-\op{hyp}(g).
$$

By using the explicit formula for $\left|\yy\right|$ of \cite[Thm. 2.2]{lmnx} together with Theorem \ref{formula} we get an explicit formula for $\op{sd}(g)$ as a polynomial in $q$ with integer coefficients.

\begin{theorem}\label{sd}
\begin{multline*}
\op{sd}(g)=\sum_{1<m|2g+2}\varphi(m)\left[A_0\left(\dfrac {2g+2}m\right)\right]_{m|q+1,\,\frac {2g+2}m \mbox{\tiny \,odd}}+\\+
\sum_{1<m|2g}\varphi(m)\left(\left[A_0\left(\dfrac {2g}m\right)\right]_{m|q+1,\,\frac {2g}m \not\equiv \frac {q+1}m\mdt2}+\left[A_2\left(\frac{2g}m\right)\right]_{m|q-1,\,\frac {q-1}m \mbox{\tiny\, odd}}\right).\end{multline*}
\end{theorem}

\begin{corollary}
If $g$ is odd and $q\equiv1\md4$ there are no self-dual hyperelliptic curves of genus $g$ defined over $k$.  
\end{corollary}

\begin{proof}
None of the conditions appearing in the formula for $\op{sd}(g)$ are satisfied. If $m|(2g+2)$ and $(2g+2)/m$ is odd then $4|m$ and $m\nmid q+1$. If $m|2g$ and $m|q+1$, then $2g/m$ and $(q+1)/m$ have the same parity. Finally, if $m|2g$ and $m|q-1$, then necessarily $(q-1)/m$ is even. 
\end{proof}

We display in Table 2 the value of $\op{sd}(g)$ for $2\le g\le 10$. As before, the contribution of $m=2$ in the formula gives a good approximation to $\op{sd}(g)$.

\begin{corollary}
$$\as{1.2}
\begin{array}{ll}
\op{sd}(g)=q^g-2q^{g-2}+2q^{g-4}+O(q^{g-6}),&\qquad \mbox{ if $g$ even, }g\ge 10,\ 4|q-1, \\
\op{sd}(g)=q^g-2q^{g-2}+2q^{g-3}+O(q^{g-6}),&\qquad \mbox{ if $g$ even, }g\ge 10,\ 4|q+1, \\
\op{sd}(g)=2q^{g-1}-2q^{g-2}+O(q^{g-5}),&\qquad \mbox{ if $g$ odd,  }g\ge 9,\ 4|q+1.
\end{array}
$$
\end{corollary}

\begin{proof}
The contribution of $m=2$ in the general formula is:
$$
[A_0(g+1)]_{g\mbox{\tiny\, even}}+[A_0(g)]_{g\not\equiv(q+1)/2\mdt2}+[A_2(g)]_{4|q+1}.
$$
The result follows by applying Lemma \ref{atilla}, having in mind that the highest power of $q$ arising from the other terms is $q^{\lceil(g-2)/2\rceil}$, corresponding to $m=4$.
\end{proof}

\begin{center}
\begin{table}
\caption{\small Number of self-dual hyperelliptic curves of genus $g$ up to $k$-isomorphism. For $g$ odd we assume that $q\equiv3\md4$. For $g$ odd and $q\equiv1\md4$ it is $\op{sd}(g)=0$}\as{1.2}
\begin{tabular}{|c|l|}\hline
\mbox{\footnotesize$g$}&\qquad\qquad\qquad\qquad\qquad\qquad\qquad\mbox{\footnotesize$\op{sd}(g)$}\\\hline
\mbox{\footnotesize$2$}&\mbox{\footnotesize$q^2-2+[2]_{3|q+1}+[2]_{q\equiv5,\,7\mdt8}$}\\\hline
\mbox{\footnotesize$3$}&\mbox{\footnotesize$2q^2-2q+[2]_{p>3}+[4]_{8|q+1}$}\\\hline
\mbox{\footnotesize$4$}&\mbox{\footnotesize$
q^4-2q^2+2+2[q-1]_{4|q+1}+2[q-1]_{q\equiv3,\,5\mdt8}+[4]_{5|q+1}+[4]_{q\equiv9,\,15\mdt{16}}$}\\\hline
\mbox{\footnotesize$5$}&\mbox{\footnotesize$2q^4-2q^3+2q^2-2q+[4]_{3|q+1}+[4]_{q\equiv\pm1\mdt5}$}\\\hline
\mbox{\footnotesize$6$}&\negmedspace\negmedspace\mbox{\footnotesize$\begin{array}{l}q^6-2q^4+2q^2-2+2[q^3-q^2]_{4|q+1}+2[q^2-q+1]_{8|q-5}+2[q^2-q-1]_{8|q+1}+\\\qquad\qquad\qquad\qquad+2[q-1]_{q\equiv5,\,7\mdt{12}}+[6]_{7|q+1}+[4]_{q\equiv13,\,23\mdt{24}}\end{array}$}\\\hline
\mbox{\footnotesize$7$}&\mbox{\footnotesize$
2q^6-2q^5+2q^2-2q+[8]_{16|q+1}+[6]_{q\equiv\pm1\mdt7}$}\\\hline
\mbox{\footnotesize$8$}&\negmedspace\negmedspace\mbox{\footnotesize$\begin{array}{l}
q^8-2q^6+2q^4-2q^2+2+2[q^5-q^4+q-1]_{4|q+1}+2[q^3-q^2+q-1]_{8|q-5}+\\\qquad+2[q^3-q^2-q+1]_{8|q-3}+2[q^2-q-1]_{3|q+1}+4[q-1]_{q\equiv7,\,9\mdt{16}}+\\
\qquad\qquad\qquad\qquad +[6]_{9|q+1}+[8]_{q\equiv17,\,31\mdt{32}}\end{array}$}\\\hline
\mbox{\footnotesize$9$}&\mbox{\footnotesize$
2q^8-2q^7+4q^4-4q^3+2-2[q^2-q-1]_{p=3}+[4]_{3|q-1}+[8]_{5|q+1}+[6]_{q\equiv\pm1\mdt9}$}\\\hline
\mbox{\footnotesize$10$}&\negmedspace\negmedspace\mbox{\footnotesize$\begin{array}{l}
q^{10}-2q^8+2q^6-2q^4+2q^2-2+2[q^7-q^6+q^3-q^2]_{4|q+1}+\\\qquad\qquad+2[q^4-q^3+q^2-q+1]_{8|q-5}+2[q^4-q^3-q^2+q+1]_{8|q+1}+\\\qquad\qquad\qquad\qquad+4[q-1]_{q\equiv9,\,11\mdt{20}}+[10]_{11|q+1}+[8]_{q\equiv21,\,39\mdt{40}}\end{array}$}\\\hline
\end{tabular}
\end{table}
\end{center}

\be

\noindent Departament de Matem\`atiques\\
  Universitat Aut\`onoma de Barcelona\\ Edifici C\\ 08193 Bellaterra,
  Barcelona, Spain\\
nart@mat.uab.cat

\end{document}